\documentclass{article}

\usepackage{amsmath}
\usepackage{amssymb}
\usepackage{epsfig}

\newcommand{\x}{\textbf{x}}
\renewcommand{\v}{\textbf{v}}
\renewcommand{\k}{\textbf{k}}
\renewcommand{\d}{\partial}
\newcommand{\R}{\mathbb{R}}
\renewcommand{\u}{\textbf{u}}
\newcommand{\eps}{\varepsilon}

\newcommand{\Qt}{\tilde{Q}}

\def \Ghat {\widehat{G}}

\begin{document}

\title{Conservative Deterministic Spectral Boltzmann Solver
near the grazing collisions limit}


\author{Jeffrey R. Haack \thanks{Department of Mathematics, The University of Texas at Austin, 2515 Speedway, Stop C1200 Austin, Texas 78712}  \and Irene M. Gamba \thanks{Department of Mathematics, The University of Texas at Austin, 2515 Speedway, Stop C1200 Austin, Texas 78712 and ICES, The University of Texas at Austin, 201 E. 24th St., Stop C0200, Austin, TX 78712} 
}

\maketitle

\begin{abstract}
We present new results building on the conservative deterministic spectral method for the space homogeneous Boltzmann equation developed by Gamba and Tharkabhushaman. This approach is a two-step process that acts on the weak form of the Boltzmann equation, and uses the machinery of the Fourier transform to reformulate the collisional integral into a weighted convolution in Fourier space. A constrained optimization problem is solved to preserve the mass, momentum, and energy of the resulting distribution. 

Within this framework we have extended the formulation to the case of more general case of collision operators with anisotropic scattering mechanisms, which requires a new formulation of the convolution weights. We also derive the grazing collisions limit for the method, and show that it is consistent with the Fokker-Planck-Landau equations as the grazing collisions parameter goes to zero.
\end{abstract}


\section{Introduction} \label{sec:intro}

There are many difficulties associated with numerically solving the
Boltzmann equation, most notably the dimensionality of the problem and
the conservation of the collision invariants. For physically relevant
three dimensional applications the distribution function is seven
dimensional and the velocity domain is unbounded. In addition, the
collision operator is nonlinear and requires evaluation of a five
dimensional integral at each point in phase space. The collision
operator also locally conserves mass, momentum, and energy, and any
approximation must maintains this property to ensure that macroscopic
quantities evolve correctly.

Spectral methods are a deterministic approach that compute the
collision operator to high accuracy by exploiting its Fourier
structure. These methods grew from the analytical works of Bobylev \cite{Bob88} developed for the Boltzmann equation with 
 Maxwell type potential interactions and integrable angular cross section, where the corresponding Fourier transformed equation has a 
closed form. Spectral approximations for these type of models where first proposed by Pareschi and Perthame \cite{ParPer96}.
Later Pareschi and Russo \cite{ParRus00} applied this work to variable hard potentials by periodizing the problem and its solution and implementing spectral collocation methods.

These methods require $O(N^{2d})$ operations per evaluation of the collision operator, where $N$ is the
total number of velocity grid points in each dimension. While convolutions 
can generally be computed in $O(N^d \log N)$ operations, the
presence of the convolution weights requires the full $O(N^{2d})$ computation of the convolution,  except for a few special cases, e.g.,
the Fokker-Planck-Landau collision operator \cite{ParRusTos00, MouPar06}. Spectral methods provide many advantages over Direct Simulation Monte Carlo Methods (DSMC) because they are more suited to time dependent problems, low Mach number flows, high mean velocity flows, and flows that are away from equilibrium. In addition, deterministic methods avoid the statistical fluctuations that are typical of particle based methods.

Inspired by the work of Ibragimov and Rjasanow \cite{IbrRja02}, Gamba
and Tharkabhushanam \cite{GamTha09, GamTha10} observed that the
Fourier transformed collision operator takes a simple form of a weighted
convolution and  developed a spectral method based
on the weak form of the Boltzmann equation that provides a general
framework for computing both elastic and inelastic collisions. Macroscopic conservation is enforced by solving a numerical constrained optimization problem that finds the closest distribution function in $L_2$ to the output of the collision term that conserves the macroscopic quantities.  

This paper presents an extension of this method to elastic collisional models that are anisotropic with respect to the scattering direction. When the angular cross section of the collision kernel has a non-integrable singularity at zero, under certain conditions one can formally calculate the Landau asymptotics \cite{Landau37, Rosenbluth57} to obtain the well known Fokker-Planck-Landau equation. In particular, this can apply to the well known Rutherford Coulombic potential \cite{Ruther11}, which is of interest for plasma physics applications. These facts combined with our spectral constrained optimization computational schemes allow us to investigate a general class of soft potential collisional kernels with  anisotropic singular angular sections as well as the grazing collision transition regime between the classical Boltzmann equation to the Fokker-Planck-Landau equation. 

\section{The space homogeneous Boltzmann equation} \label{sec:BTE}
The space homogeneous Boltzmann equation is given by the initial value problem
\begin{equation}\label{BTE}
\frac{\d}{\d t} f(\v,t) = \frac{1}{Kn} Q(f,f),
\end{equation}
with
\begin{align*}
\qquad \v \in \R^d, \qquad f(v,0) = f_0(\v) \\
\end{align*}
where $f(\v,t)$ is a probability density distribution in $\v$-space and  $f_0$ is assumed to be locally integrable with respect to $\v$. The dimensionless  parameter $Kn > 0$ is the scaled Knudsen number, which is defined as the ratio between the mean free path between collisions and a reference macroscopic length scale.
 
The collision operator $Q(f,f)$ is a bilinear integral form in $(\v,t)$ given by
\begin{equation}\label{Q_general}
Q(f,f)(\v,t) = \int_{\v_\ast \in \R^d} \int_{\sigma \in S^{d-1}} B(|\v - \v_\ast|,\cos \theta) (f(\v_\ast ')f(\v') - f(\v_\ast)f(\v)) d\sigma d\v_\ast,
\end{equation}
where the velocities $\v', \v_\ast'$ are determined through a given collision rule \eqref{velocity_interact}, depending on $\v, \v_\ast$. The positive term of  the integral in \eqref{Q_general} evaluates  $f$ in the pre-collisional velocities that can result in a post-collisional velocity the direction $\v$.  The collision kernel $B(|\v - \v_\ast|,\cos \theta)$ is a given non-negative function depending on the size of the relative velocity $\u := \v - \v_\ast$ and $\cos \theta = \frac{\u \cdot \sigma}{|\u|}$, where $\sigma$ in the $n-1$ dimensional sphere $S^{n-1}$ is referred ro as the scattering direction which also coincides with the direction of the post-collisional elastic relative velocity.

For the following we will use the elastic (or reversible) interaction law
\begin{align}
&\v' = \v + \frac{1}{2}(|u|\sigma - \u), \qquad \v_\ast ' = \v_\ast - \frac{1}{2}(|u|\sigma - \u) \label{velocity_interact} \\
&B(|u|,\cos \theta) = |u|^\lambda b(\cos \theta)  \ .\qquad \qquad  \nonumber
\end{align}
The angular cross section function $b(\cos \theta)$ may or may no be integrable with respect to $\sigma$ on $S^{d-1}$. When integrability holds is referred to as the Grad cut-off assumption on the angular cross section. 

The parameter $\lambda$ regulates the collision frequency as a function of the relative speed $|\u|$. This corresponds to the interparticle potentials used in the derivation of the collisional kernel and are referred to as variable hard potentials (VHP) for $0 < \lambda < 1$, hard spheres (HS) for $\lambda = 1$, Maxwell molecules (MM) for $\lambda = 0$, and variable soft potentials (VSP) for $-3 < \lambda < 0$. The $\lambda = -3$ case corresponds to a Coulombic interaction potential between particles. If $b(\cos \theta)$ is independent of $\sigma$ we call the interactions isotropic, e.g, in the case of hard spheres in three dimensions.

\subsection{Spectral formulation} \label{sec:spectral_cont}
The key step our formulation of the spectral numerical method is the use of the weak form of the Boltzmann collision operator. For a suitably smooth test function $\phi(\v)$ the weak form of the collision integral is given by
\begin{equation} \label{collision_weakform}
\int_{\R^d} Q(f,f) \phi(\v) d\v = \int_{\R^d \times \R^d \times S^{d-1}} f(\v)f(\v_\ast) B(|\u|,\cos \theta) (\phi(\v') - \phi(\v)) d\sigma d\v_\ast d\v
\end{equation}
If one chooses 
\begin{equation*}
\phi(\v) = e^{-i \zeta \cdot \v} / (\sqrt{2\pi})^d,
\end{equation*}
then \eqref{collision_weakform} is the Fourier transform of the collision integral with Fourier variable $\zeta$:
\begin{align}\label{FourierQ}
\widehat{Q}(\zeta) &= \frac{1}{(\sqrt{2\pi})^d} \int_{\R^d} Q(f,f) e^{-i \zeta \cdot \v} d\v \nonumber \\
&= \int_{\R^d \times \R^d \times S^{d-1}} f(\v)f(\v_\ast) \frac{B(|\u|,\cos \theta)}{(\sqrt{2\pi})^d} (e^{-i \zeta \cdot \v'} - e^{-i \zeta \cdot \v}) d\sigma d\v_\ast d\v \nonumber\\
&= \int_{\R^d} G(\u,\zeta) \mathcal{F}[f(\v)f(\v-\u)](\zeta) d\u,
\end{align}

where $\widehat{[\cdot]} = \mathcal{F}(\cdot)$ denotes the Fourier transform and 
\begin{equation} \label{G_eqn}
G(\u,\zeta) = |u|^\lambda \int_{S^{d-1}} b(\cos \theta) \left(e^{-i\frac{\beta}{2} \zeta \cdot |u|\sigma}e^{i\frac{\beta}{2} \zeta \cdot \u} - 1\right) d\sigma
\end{equation}
Further simplification can be made by writing the Fourier transform inside the integral as a convolution of Fourier transforms:
\begin{align} \label{Cont_spectral}
\widehat{Q}(\zeta) 
&= \int_{\R^d} \widehat{G}(\xi,\zeta) \hat{f}(\zeta - \xi) \hat{f}(\xi) d\xi,
\end{align}

where the convolution weights $\widehat{G}(\xi,\zeta)$ are given by
\begin{align} \label{Ghat_eqn}
\widehat{G}(\xi,\zeta) &= \frac{1}{(\sqrt{2\pi})^d}  \int_{\R^d} G(\u,\zeta) e^{-i \xi \cdot u} d\u \\
&=\frac{1}{(\sqrt{2\pi})^d}  \int_{\R^d} |u|^\lambda e^{-i\xi\cdot\u} \int_{S^{d-1}} b(\cos \theta)  \left(e^{-i\frac{\beta}{2} \zeta \cdot |u|\sigma}e^{i\frac{\beta}{2} \zeta \cdot \u} - 1\right) d\sigma d\u \nonumber
\end{align}
These convolution weights can be precomputed once to high accuracy and stored for future use. For many collisional models, such as isotropic collisions, the complexity of the integrals in the weight functions can be reduced dramatically through analytical techniques. However in this paper we make no assumption on the isotropy of $b$ and derive a more general formula.

We begin with $G(\u,\zeta)$ and define a spherical coordinate system for $\sigma$ with a pole in the direction of $\u$ and obtain
\begin{align}
G(\textbf{u}, \zeta) &=
2\pi|u|^\lambda \int_0^{\pi} b(\cos \theta) \sin\theta \left( e^{i\frac\beta{2} (1-\cos\theta)\zeta \cdot \textbf{u}}J_0\left(\frac{\beta|u|\sin\theta|\zeta^\perp|}{2}\right)  - 1\right)  d\theta, \label{SphG}
\end{align}
where $J_0$ is the Bessel function of the first kind, $j,k$ are unit length basis vectors that are mutually orthogonal to $\u$, and $\zeta^\perp = \zeta - (\zeta \cdot \textbf{u}/|u|) \textbf{u}/|u|$. Note that for the isotropic case $b(\cos \theta)$ drops out and one can take $\zeta$ as the polar direction for $\sigma$, resulting in an explicit expression involving a sinc function.

$\widehat{G}$ is the Fourier transform of $G$, however we must take this transform on a {\em ball centered at 0} rather than simply taking the Fast Fourier Transform (FFT) \cite{FFTW} of $G$, which ensures that the weights are real-valued.

The convolution weights $\Ghat(\zeta,\xi)$ in 3 dimensions are computed as follows
\begin{align*}
\Ghat(\zeta,\xi) &= 2\pi \int_{B_L(0)} |u|^\lambda e^{-i\xi \cdot \u} \int_0^\pi b(\cos\theta)\sin\theta \left[e^{\frac{i\zeta}{2} \cdot \u(1-\cos\theta)} J_0\left(\frac12 |u||\zeta^\perp|\sin\theta\right) - 1\right] d\theta d\u \\
&=2\pi \int_0^L  \int_{S^2} r^{\lambda+2} \int_0^\pi b(\cos\theta)\sin\theta \left[e^{-ir(\xi - \frac{\zeta}{2}(1-\cos\theta) )\cdot \eta}J_0\left(\frac12 r|\zeta^\perp|\sin\theta\right)  -  e^{-ir\xi \cdot \eta} \right] d\theta d\eta dr.
\end{align*}

We now take $\zeta$ to be the polar direction for the spherical integration of $\eta$ and use that $\widehat{G}$ is real-valued to obtain

\begin{align}
\Ghat(\zeta,\xi) &=
4\pi^2 \int_0^L  r^{\lambda+2} \int_0^\pi \int_0^\pi b(\cos\theta)\sin\theta\sin\phi J_0\left(r\left|\xi - \frac{\xi\cdot\zeta}{|\zeta|^2}\zeta\right|\sin\phi\right) \times \nonumber \\
 & \left[\cos\left(r(\xi - \frac{\zeta}{2}(1-\cos\theta) )\cdot \frac{\zeta}{|\zeta|}\cos\phi\right) J_0\left(\frac12 r|\zeta|\sin\phi\sin\theta\right) - \cos\left(r\xi \cdot \frac{\zeta}{|\zeta|}\cos\phi\right) \right] d\theta d\phi dr. \label{GHat_calc}
\end{align}

This requires a three dimensional integral for the $N^6$ pairs $(\zeta, \xi)$, which is two orders of magnitude more than the isotropic case, but as before this weight is precomputed only once and used in any subsequent computations with the collisional model.


\section{The Conservative Numerical Method} \label{sec:numerics_setup}

\subsection{Velocity space discretization} \label{sec:discretization}

In order to compute the Boltzmann equation we must work on a bounded velocity space, rather than all of $\R^d$. However typical distributions are supported on the entire domain, for example the Maxwellian equilibrium distribution. Even if one begins with a compactly supported initial distribution, each evaluation of the collision operator spreads the support by a factor of $\sqrt{2}$, thus we must use a working definition of an {\em effective support} of size  $R$ for the distribution function. Bobylev and Rjasanow \cite{BobRja99} suggested using the temperature of the distribution function, which typically decreases as $\text{exp}(-|v|^2 / 2T)$ for large $|v|$, and used a rough estimate of $R \approx 2\sqrt{2}T$ to determine the cutoff. Thus, we assume that the distribution function is negligible outside of a ball 
\begin{equation} \label{Ball_domain_v}
B_{R_x}(\textbf{V}(\x)) = \{\v \in \R^d : |\v - \textbf{V(\x)}| \le R_x \},
\end{equation}
where $\textbf{V}(\x)$ is the local flow velocity which depends in the spatial variable $\x$. For ease of notation in the following we will work with a ball centered at $0$ and choose a length $R$ large enough that $B_{R_x}(\textbf{V}(\x)) \subset B_R(0)$ for all $\x$.

With this assumed support for the distribution $f$, the integrals in \eqref{Cont_spectral} will only be nonzero for $\u \in B_{2R}(0)$. Therefore, we set $L=2R$ and define the cube
\begin{equation} \label{Cube_domain_v}
C_L = \{ \v \in \R^d : |v_j| \le L,\,\, j = 1,\dots,d\}
\end{equation}
to be the domain of computation. With this domain the comptuation of the weight function integral \eqref{GHat_calc} is cut off at $r_0=L$.

Let $N \in \mathbb{N}$ be the number of points in velocity space in each dimension. Then we establish a uniform velocity mesh with $\Delta v = \frac{2L}{N}$ and due to the formulation of the discrete Fourier transform the corresponding uniform Fourier space mesh size is given by $\Delta \zeta = \frac{\pi}{L}$. 
%


\subsection{Collision step discretization} \label{sec:collision_disc}
%

To simplify notation we will use one index to denote multidimensional sums with respect to an index vector $\textbf{m}$
\begin{equation*}
\sum_{\textbf{m}=0}^{N-1} = \sum_{m_1,\dots,m_d = 0}^{N-1}.
\end{equation*}

To compute $\widehat{Q}(\zeta_\k)$, we first compute the Fourier transform integral giving $\hat{f}(\zeta_k)$ via the FFT. The weighted convolution integral is approximated using the trapezoidal rule
\begin{equation}
\widehat{Q}(\zeta_\k)= \sum_{\textbf{m} = 0}^{N-1} \widehat{G}(\xi_{\textbf{m}},\zeta_\k) \hat{f}(\xi_{\textbf{m}}) \hat{f}(\zeta_\k - \xi_{\textbf{m}}) \omega_{\textbf{m}},
\end{equation}
where $\omega_\textbf{m}$ is the quadrature weight and we set $\hat{f}(\zeta_\k - \xi_{\textbf{m}}) = 0$ if $(\zeta_\k - \xi_{\textbf{m}})$ is outside of the domain of integration. We then use the inverse FFT on $\widehat{Q}$ to calculate the integral returning the result to velocity space. 

Note that in this formulation the distribution function is not periodized, as is done in the collocation approach of Pareschi and Russo \cite{ParRus00}. This is reflected in the omission of Fourier terms outside of the Fourier domain. All integrals are computed directly only using the FFT as a tool for faster computation and the convolution integral is accurate to at least the order of the quadrature. The calculations below use the trapezoid rule, but in principle Simpson's rule or some other uniform grid quadrature can be used. However, it is known that the trapezoid rule is spectrally accurate for periodic functions on periodic domains (which is the basis of spectral accuracy for the FFT), and the same arguments can apply to functions with sufficient decay at the integration boundaries \cite{Atkinson}. These accuracy considerations will be investigated in future work. The overall cost of this step is $O(N^{2d})$. 

\subsection{Discrete conservation enforcement} \label{sec:conservation}
This implementation of the collision mechanism does not conserve all of the quantities of the collision operator. To correct this, we formulate these conservation properties as a Lagrange multiplier problem. Depending on the type of collisions we can change this constraint set (for example, inelastic collisions do not preserve energy), but we will focus on the case of elastic collisions, which preserve mass, momentum, and energy. 

Let $M = N^d$ be the total number of grid points, let $\tilde{\textbf{Q}} = (\Qt_1, \dots, \Qt_M) ^T$ be the result of the spectral formulation from the previous section, written in vector form, and let $\omega_j$ be the quadrature weights over the domain in this ordering. Define the integration matrix
\begin{equation*}
\textbf{C}_{5\times M} = \left(\begin{array}{c} \omega_j \\ v_j^i \omega_j \\ |\v_j|^2 \omega_j \end{array} \right),
\end{equation*}
where $v^i,\, i=1,2,3$ refers to the $i$th component of the velocity vector. Using this notation, the conservation method can be written as a constrained optimization problem. 

\begin{equation} 
\text{Find } \textbf{Q} = (Q_1,\dots,Q_M)^T \text{ that minimizes } \frac12 \|\tilde{\textbf{Q}} - \textbf{Q}\|_2^2 \text{ such that } \textbf{C} \textbf{Q} = \textbf{0}
\end{equation}
Formulating this as a Lagrange multiplier problem, we define 
\begin{equation}
L(\textbf{Q},\gamma) = \sum_{j=1}^M (\Qt_j - Q_j)^2 - \gamma^T\textbf{C}\textbf{Q}
\end{equation}
The solution is given by
\begin{align}
\textbf{Q} &= \tilde{\textbf{Q}} + \textbf{C}(\textbf{C}\textbf{C}^T)^{-1} \textbf{C} \tilde{\textbf{Q}} \nonumber \\
&:= \textbf{P}_N \tilde{\textbf{Q}}
\end{align}

Overall the collision step in semi-discrete form is given by
\begin{equation}
\frac{\d \textbf{f}}{\d t} = \textbf{P}_N \tilde{\textbf{Q}}
\end{equation}

The overall cost of the conservation portion of the algorithm is a $O(N^d)$ matrix-vector multiply, significantly less than the computation of the weighted convolution.

\subsection{Grazing collisions limit}
The Fokker-Planck-Landau equation is used to describe binary collisions occurring in a plasma, for example, and can be shown to be an approximation of the Boltzmann equation in the case where the dominant collision mechanism is that of grazing collisions, i.e., collisions that only result in very small deflections of particle trajectories, as is the case for Coulomb potentials with Rutherford scattering \cite{Ruther11} between charged particles. The equation is given by
\begin{equation}\label{FPL}
Q_L (f,f) = \nabla_v \cdot \left( \int_{\mathbb{R}^3} |u|^{\lambda+2} (I - \frac{\u\otimes \u}{|u|^2}) (f(v_\ast)\nabla_v f(v) - f(v) (\nabla_v f)(v_\ast)) d v_\ast\right) \, .
\end{equation}

Applying the same methodology used to derive the spectral method for the Boltzmann equation, the Fourier transform of the Fokker-Planck-Landau equation takes the form \eqref{FourierQ}, where the  weight function $G(\u,\zeta)$ in \eqref{G_eqn} is now given by
\begin{equation} \label{GFPL} 
G(\u,\zeta) = 4 |u|^\lambda (i(\u\cdot\zeta) - \frac14 |u|^2|\zeta^\perp|^2) \, .
\end{equation}

To show that the spectral method for Boltzmann operator is consistent with this form of Fokker-Planck-Landau operator, we must take the grazing collisions limit, which requires that the angular scattering function is consistent with the singular rates of Rutherford scattering. Indeed, it is enough to assume that the collision kernel satisfies the following.

Let $\eps > 0$ be the small parameter associated with the grazing collisions limit. A family of kernels $b_\eps(\theta)$ represents grazing collisions if:
\begin{itemize}
 \item $\lim_{\eps \to 0}  \Lambda_{\eps} = 2\pi \int_0^{\pi} b(\cos \theta , \eps) (1-\cos\theta)\, 
\sin\theta d \theta =\Lambda_0 <\infty$  

 \item $\lim_{\eps \to 0} = 2\pi \int_0^{\pi} b(\cos \theta , \eps) (1-\cos\theta)^{2+k}\, 
\sin\theta d \theta \to 0 \qquad \mathrm{for}\ \ k\geq 0 \, .$
\end{itemize}

In particular, the angular part of the classical Rutherford scattering for $\theta$ near 0 corresponds  to a  
family $b_{\eps}(\cos\theta)$ given by
\begin{equation}
b(\cos \theta, \eps) \sin\theta =   \frac{8 \eps }{\pi \theta^4} 1_{\theta \ge \eps} \label{beps} \, , 
\end{equation}
 which satisfies the above conditions for $\Lambda_0=8$. This is  the potential we take for the following calculation and numerical simulation of the grazing collision limit for the Boltzmann equation for Coulombic interactions in three dimensions, i.e., $\lambda=-3$.

For the grazing collisions limit, take Taylor expansion of the exponential term in $G$ \eqref{SphG}
\begin{align*} 
&e^{i\frac\beta{2}((1-\cos\theta)\zeta\cdot\textbf{u} + |u|\zeta\cdot j \sin\theta \sin\phi + |u|\zeta\cdot k\sin\theta\cos\phi)} - 1 = i\frac\beta{2}((1-\cos\theta)\zeta\cdot\textbf{u} + |u|\zeta\cdot j \sin\theta \sin\phi + |u|\zeta\cdot k\sin\theta\cos\phi) \\
&-\frac{\beta^2}8 ((1-\cos\theta)\zeta\cdot\textbf{u} + |u|\zeta\cdot j \sin\theta \sin\phi + |u|\zeta\cdot k\sin\theta\cos\phi)^2 
+ O(|\zeta|^3|u|^3\theta^3)
\end{align*}
We can toss the $\theta^3$ terms because as $\eps \to 0$ 
\[\int_0^\pi b_\eps(\cos \theta)\sin \theta\, \theta^3 d\theta \to 0 \]

Taking the $\phi$ integral of the first order term gives
\[i\pi\beta(1-\cos\theta)\zeta \cdot \textbf{u}. \]
For the second order term, taking the $\phi$ integral results in, after removing integrals that evaluate to zero by symmetry,
\[-\frac{\beta^2}{8} \int_0^{2\pi} ((1-\cos\theta)^2(\zeta\cdot\textbf{u})^2 + |u|^2(\zeta\cdot j)^2 \sin^2\theta \sin^2 \phi + |u|^2(\zeta \cdot k)^2\sin^2\theta \cos^2\phi)\]
We can neglect the $(1-\cos\theta)^2$ term as it is $O(\theta^4)$. Thus we obtain the small $\eps$ approximation for $G$
\begin{align*}
G(\textbf{u},\zeta) &= 
 \pi|u|^\lambda\int_0^\pi b_\eps(\cos \theta) \sin \theta \left(i\beta(1-\cos\theta)\zeta\cdot\textbf{u} - \frac{|u|^2|\zeta^\perp|^2\beta^2}{8}\sin^2 \theta\right)
\end{align*}
From the earlier definition of the collision kernel, we have that
\[ 2\pi\sin\theta(1-\cos\theta)b_\eps(\cos\theta) \to \Lambda_0\, \delta_{\theta=0}\]
Thus
\begin{equation}
G(\textbf{u},\zeta) = \frac{\Lambda_0}{2} |u|^\lambda \left(i\beta(\zeta \cdot \textbf{u}) - \frac{|u|^2|\zeta^\perp|^2\beta^2}{4}\right)
\end{equation}
When $\Lambda_0=8$, this gives the same $G(\u,\zeta)$ as found above by applying the spectral technique to the Fokker-Planck-Landau equation.

\subsection{Computing $\Ghat$ for singular scattering kernels}

Numerically calculating the weights $\hat{G}$ to high accuracy can be difficult for singular scattering kernels, due to the precise nature of the cancellation at the left endpoint of the integral. Using $b_\eps(\cos \theta)$ from \eqref{beps}, the $\theta$ integral in \eqref{GHat_calc} is given by
\begin{equation}
\int_\eps^\pi \frac{8\eps}{\pi \theta^4} (\cos(c_1 (1-\cos\theta) - c_3)  J_0(c_2 \sin \theta) - cos(c_3)) d\theta,
\end{equation}
where $c_1, c_2, c_3$ depend on the current values of $\phi, r, \zeta, \xi$ following from the full formulation of $\Ghat$. When $\eps << 1$ the bulk of the integration mass occurs near the left endpoint of the $\theta$ interval, however this presents a challenge for a numerical quadrature package to compute. For $\theta << 1$ there is a subtraction of two nearly equal numbers, which causes floating point errors, as well as a multiplication of two numbers of very different magnitudes ($\eps$ and $\theta^{-4}$). To alleviate this, we split the integration interval into two pieces, and use the first term of the Taylor expansion of the integrand for $\theta << 1$:
\begin{equation}
\int_\eps^{\sqrt{\eps}} \frac{8\eps}{\pi \theta^2} \left(-\frac{c_2^2}{4}\cos(c_3) + \frac{c_1}{2}\sin(c_3)\right)d\theta + 
\int_{\sqrt{\eps}}^\pi \frac{8\eps}{\pi \theta^4} (\cos(c_1 (1-\cos\theta) - c_3)  J_0(c_2 \sin \theta) - \cos(c_3)) d\theta.
\end{equation}
These integrals are computed using the GNU Scientific Laboratory integration routines \cite{gsl}. We use \verb+cquad+ to compute the first $\theta$ integral, which has proved to be most stable choice for this near-singular integrand. The adaptive Gauss-Konrod quadrature \verb+qag+ is used for all other integrals used in computing the weights $\Ghat$.
\section{Numerical results}
\begin{figure}[!htbp] \label{GrazingSlice}
\includegraphics[width=.6\linewidth]{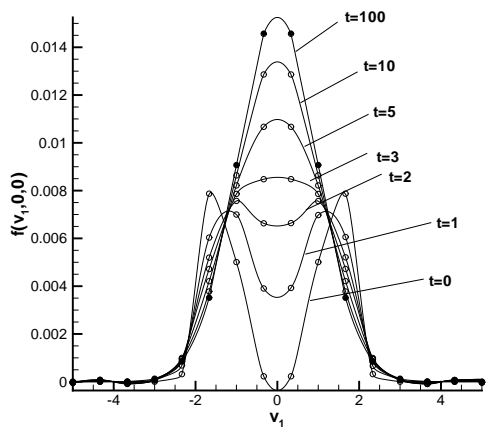}
\caption{Slice of distribution function at different times. The solid dots are the values of the Maxwellian distribution associated with the initial data. The solid lines are a spline reconstruction based on the grid values (hollow circles).}
\end{figure}

To illustrate that this method captures the correct behavior for grazing collisions, we take a Coulombic potential ($\lambda = -3$) and set $\eps = 10^{-4}$. Similar to what was done in \cite{ParRusTos00} for the Landau equation, we take the axially symmetric initial condition
\[ f(\v,0) = 0.01  \textrm{exp}\left(-10 \left(\frac{|\v| - 1.5}{1.5}\right)^2\right). \]
We take a domain size of $L = 5$ and compute to time $t = 100$ with a timestep of $0.001$. The results are shown in Figure \ref{GrazingSlice}.
\section{Conclusions and future work}
We have derived the precomputed weight formulas for the more generalized anisotropic collisional models within the Boltzmann equation. We also showed that the spectral method for the Boltzmann equation is consistent with the limiting Fokker-Planck-Landau equation under suitable assumptions on the scattering kernel. One other important note is that this method could be a good candidate for collisional models where the collision mechanism is unknown and only experimentally determined. In addition, as the Landau equation is used to model collisions of charged particles in plasma we will seek to add field effects to the space inhomogeneous Boltzmann equation, resulting in the Boltzmann-Poisson or Boltzmann-Maxwell systems. The inhomogeneous method uses operator splitting between the collision and the transport terms, so in principle one can use whatever transport solver one wants for the spatial terms in the equation. 

\section{Acknowledgments}
This work has been supported by the NSF under grant number DMS-0636586.





\bibliographystyle{siam}   

\bibliography{Boltz}


\end{document}